\documentclass[12pt,reqno,twoside]{amsart}
\usepackage{amsmath}
\usepackage{amssymb}
\usepackage{graphicx}

\newtheorem {Theorem}    {Theorem}[section]

\newtheorem {Conjecture} [Theorem]    {Conjecture}

\theoremstyle{definition}
\newtheorem {Definition} [Theorem]    {Definition}

\newtheorem {Question}   [Theorem]    {Question}
\newcounter{AbcT}

\numberwithin{equation}{section}

\newcommand {\equ}[1]     {\eqref{#1}}

\newcommand {\cond}   {\:|\:}

\newcommand {\N} {{\mathbb N}} 
\newcommand {\Q} {{\mathbb Q}} 
\newcommand {\R} {{\mathbb R}}
 
\newcommand {\T} {{\mathbb T}} 
\renewcommand {\H} {{\mathbb H}}
\newcommand {\Z} {{\mathbb Z}}

\renewcommand{\liminf}{\varliminf}

\DeclareMathOperator{\supp}{supp}

\DeclareMathOperator{\image}{Image}

\DeclareMathOperator{\SL}{SL}

\DeclareMathOperator{\GL}{GL}

\DeclareMathOperator{\PGL}{PGL}

\newcommand {\Borel}      {\mathcal B}

\newcommand {\IGNORE}[1]  {}

\newcommand {\absolute}[1] {\left| {#1} \right|}
\newcommand {\norm}[1] {\left\| {#1} \right\|}
\newcommand{\prefixitemswith}[1] {\renewcommand{\theenumi}{#1\arabic{enumi}}}

\newcommand{\twobytwo}[4] {\begin{pmatrix} {#1}&{#2}\\{#3}&{#4}\end{pmatrix}}

\newcommand{\isomof}[1] {\operatorname{Isom} (#1)}

\newcommand {\haaron}[1] {{\mathcal{H}_{\uppercase{#1}}}}

\renewcommand {\image} {\operatorname{Image}}
\begin{document}

\title{Rigidity of multiparameter actions}

\author{Elon Lindenstrauss}

\address{Courant Institute of Mathematical Sciences, New York University, New York, NY 10012}

\date{version of \today}
\email{elonbl@member.ams.org}

\maketitle

\setcounter{section}{-1}
\section {Prologue}
In this survey I would like to expose recent developments and applications of the study of the rigidity properties of natural algebraic actions of multidimensional abelian groups. This study was initiated by Hillel Furstenberg in his landmark paper \cite{Furstenberg-disjointness-1967}.

The survey is an outgrowth of notes I gave to a series of lectures in the  II Workshop on Dynamics and Randomness in December 2002; a related talk was given also in the Conference on Probability in Mathematics in honor of Hillel Furstenberg. It is not intended to be an exhaustive survey, and is skewed towards the  problems I have worked on and which I know best. In particular, I make no mention of $ \Z ^ d$ actions on zero dimensional compact groups such as Ledrappier's example of a 2-mixing but not 3-mixing $ \Z ^2$-action \cite [Intr.]{ Schmidt-book} --- a deficiency that is partially mitigated by the existence of several good surveys, for example \cite{Schmidt-PIMS-survey}.

It has been my good fortune to have been able to learn ergodic theory from Benjamin Weiss, Hillel's long term collaborator and friend, and from Hillel himself during my graduate studies in the Hebrew University. To the problems mentioned here I have been introduced by numerous discussions (as well as collaborations) with David Meiri, Shahar Mozes, Yuval Peres and Barak Weiss, all of whom were students of Hillel, as well as with Hillel himself.  It is therefore with great pleasure that I dedicate this survey to Hillel, on occasion of his retirement. After graduating, I again had the good fortune of learning from Peter Sarnak about  the quantum unique ergodicity question, which is fortuitously connected to precisely these kind of questions. Without Peter's help and encouragement, in times even his insistence, I would not have been able to discover this connection to quantum unique ergodicity which motivated much of my work on this subject.

\section{Rigidity on the one torus}

One of the simplest dynamical systems is the map
\begin{equation*}
\times_n : x \mapsto n x \bmod 1
\end{equation*}
on the unit interval, where $n$ is any natural number. In order to make
this map continuous, we think of it as a map on the 1-torus $\T = \R /
\Z$. This map is essentially the same as the shift on a one-sided
infinite sequence of $n$ symbols: the exact relationship between these
two systems is that there is a continuous surjection from the one-sided
shift on $n$ symbol to $\T$ commuting with the respective $\Z^+$
actions, and this continuous surjection is finite-to-1 everywhere (in
fact at most 2-to-1) and is 1-to-1 except at countably many points.

There are very many closed invariant sets for the map $\times_n$, and
very many invariant or even ergodic measures.
As discovered by Furstenberg, the situation changes dramatically if one
adds an additional map into the picture.
One may consider the two commuting maps $\times_n : x \mapsto nx \bmod
1$ and $\times_m : x \mapsto mx \bmod 1$; in order to really have two
different maps it is natural to assume that $n$ and $m$ are {\bf
multiplicative independent}, i.e. that $\log n / \log m \not \in \Q$.

In this case Furstenberg proved that there are very few closed
invariant sets \cite{Furstenberg-disjointness-1967}, and conjectured
that a similar statement holds for invariant measures.

\begin{Theorem} [Furstenberg \cite{Furstenberg-disjointness-1967}]
\label{Furstenberg's theorem}
Let $n , m$ be two multiplicative independent positive integers. Assume
that $X \subset \T$ is closed and invariant under $\times_n$ and
$\times_m$ (i.e. that $\times_n ( X ) \subset X$ and similarly for
$\times_m$). Then either $X = \T$ or $X$ is finite (and so necessarily
$X \subset \Q$).
\end{Theorem}

His proof is a true gem, and uses a fundamental idea, which was new at
the time, namely disjointness of dynamical systems: some systems are so
different they cannot be coupled in any nontrivial way. Like most of
what we consider here, this notion applies both in the measurable
category and in the topological category; here is the topological
version:

\begin{Definition}
Let $G$ be a reasonable semigroup (say locally compact, second
countable). Suppose $G$ acts continuously on two spaces $X$ and $Y$. A
{\bf joining} of $( X , G )$ and $( Y , G )$ is a closed subset $Z
\subset X \times Y$ invariant under the diagonal $G$ action on $X
\times Y$ whose projection to the first coordinate is $X$ and to the
second coordinate is $Y$.
The two systems $( X , G )$ and $( Y , G )$ are said to be {\bf
disjoint} if the only joining between them is the trivial product
joining $Z = X \times Y$.
\end{Definition}

In Furstenberg's proof of Theorem \ref{Furstenberg's theorem} a key
point is that a  continuous $\Z^d$-action on a compact space with many
periodic points is disjoint from any minimal $\Z^d$-action (i.e. one in
which every orbit is dense).

We now present Furstenberg's conjecture about measure rigidity, which
despite substantial progress is still open and considered one of the
most outstanding problems in ergodic theory. We recall that a measure
$\mu$ on a measurable space $( X , \Borel )$ on which a semigroup $G$
acts is said to be {\bf invariant} if for every $g \in G$ and any $A
\in \Borel$ we have that $\mu ( g^{-1} A ) = \mu ( A )$; an invariant
measure is said to be ergodic if every set $A \in \Borel$ which is
invariant under every $g \in G$ satisfies that $\mu ( A ) = 0$ or $\mu
( X \setminus A ) = 0$.

\begin{Conjecture} [Furstenberg (1967)]
Let $n , m$ be two multiplicative independent positive integers. Any
Borel measure $\mu$ on $\T$ ergodic under the action of the semigroup
generated by $\times_n , \times_m$ is either Lebesgue measure or an
atomic measure supported on finitely many rational points.
\end{Conjecture}

The first substantial partial result in this direction is due to R.
Lyons \cite{Lyons-2-and-3}. He proved the following:

\begin{Theorem}[Lyons \cite{Lyons-2-and-3}]
Let $m , n$ be two multiplicative independent positive integers, and
$\mu$ a $\times_n , \times_m$-and ergodic measure on $\T$ as in
Furstenberg's conjecture. Assume further that with respect to the $\Z^{
+ }$ action generated by the $\times_n$ map, $\mu$ has the $K$ property
(i.e. $( \T , \Borel , \mu , \times_n )$ has no zero entropy factors).
Then $\mu$ is Lebesgue measure.
\end{Theorem}

This theorem has been substantially strengthened by Rudolph, whose
theorem (at least for the case of $m , n$ relatively prime) still
presents the state-of-the-art regarding Furstenberg's conjecture:
\begin{Theorem} [Rudolph \cite{Rudolph-2-and-3}]
\label{Rudolph theorem}
Let $m , n$ be relatively prime positive integers, and $\mu$ a
$\times_n , \times_m$-ergodic measure on $\T$ as before. Assume that
the entropy of the $\Z^+$-action generated by one of these maps, say
$\times_n$, has positive entropy. Then $\mu$ is Lebesgue measure.
\end{Theorem}
Rudolph's proof used in an essential way that n and m are relatively
prime. This restriction was lifted by A. Johnson
\cite{Johnson-invariant-measures} who extended
Theorem \ref{Rudolph theorem}
to $n$ and $m$ multiplicatively independent.

In his proof, Rudolph identified and dealt with a substantial
complication, namely that even though $\mu$ is $\times n , \times
m$-ergodic it does not need to be $\times n$-ergodic. Overcoming this
requires a new idea, and if one makes this ergodicity assumption,
Rudolph's proof can be simplified. Note that by definition this
difficulty does not exist if the action of $\times n$ on $( \T , \Borel
, \mu )$ is K, hence in particular ergodic, which is a case covered by
Lyons' theorem.

An important advance was made by Katok and Spatzier
\cite{Katok-Spatzier, Katok-Spatzier-corrections}, who discovered that
Rudolph's proof can be extended to give partial information on
invariant measures in much greater generality. Here overcoming the need
for an ergodicity assumption becomes substantially harder; an extension
of Rudolph's original method of overcoming this difficulty works only
in a rather restricted class of examples (of course much more general
than that considered in Rudolph's original proof) with satisfy a
condition Katok and Spatzier called {\bf total non-symplecticity
(TNS)}.

We remark that J. Feldman has been able to give an alternative
proof of Rudolph's and Johnson's theorems using a method similar to
that of Lyons \cite{Feldman-generalization}; however this proof does
not seem very suitable for generalizations.
A third approach to Rudolph's Theorem by Host
\cite{Host-normal-numbers, Host-multidimensional} will play a key role
later on.

In this survey we emphasize recent progress in the study of measures invariant under algebraic multiparameter abelian actions, in particular progress towards the classifying of invariant measures without need for an ergodicity assumption in a much wider class of actions than the class of TNS actions for which this was previously known. It turns out that in applications to number theoretic and other problems this improvement is crucial. Also, following Host's approach mentioned above the scope of these techniques have been extended to deal with certain ``partially isometric'' {\bf one} parameter actions.

As we will explain in \S\ref{section: quantum unique ergodicity}, this study of invariant measure yields a proof of arithmetic cases of the quantum unique ergodicity conjecture \cite{Lindenstrauss-quantum}. We also mention that in \cite{Einsiedler-Katok-Lindenstrauss} (which is joint with M. Einsiedler and A. Katok)  we use a partial classification of invariant measures to derive a substantial partial result towards the long-standing Littlewood Conjecture on simultaneous Diophantine approximations (Theorem~\ref{theorem regarding Littlewood} below). In both cases measure rigidity techniques prove results which are at present unapproachable by other methods.

\section{rigidity for toral automorphism}
In this and the next section, we present two classes of actions of
multidimensional abelian groups where similar rigidity phenomena can be
found.

Much of what has been said in the first section regarding the
$\times_m$ map is also true for higher dimensional toral automorphisms,
or more generally for automorphisms of connected abelian groups. If
$\alpha$ is a single toral automorphism, and assuming the automorphism
is expansive (i.e. hyperbolic), one can find a sofic system (i.e. a
factor of a shift of finite type) which is related to the toral
automorphism in a way which is very similar to how the one-sided shift
on $n$ symbols is related to the $\times_n$ map, and so the abundance
of closed invariant sets and invariant measures on a sofic system
implies that the same holds for a hyperbolic toral automorphism.

The construction of such sofic models is a fascinating story in its own
right. A geometric construction in the two dimensional case was
discovered by Adler and Weiss \cite{Adler-Weiss}, and this theory, the
theory of Markov partitions, was developed by Sinai, Bowen and others
as an exceedingly powerful tool for the study of general Anosov and
axiom A systems. More recently, Kenyon and Vershik
\cite{Kenyon-Vershik} (following earlier work of Vershik)
constructed an algebraic coding specifically for
hyperbolic toral automorphism which better reflects the algebraic
structure
of this system; a related but different coding has been developed by K.
Schmidt \cite{Schmidt-coding}, and more careful analysis of certain
special cases was carried out by Sidorov and Vershik
\cite{Sidorov-Vershik}.

In general, automorphisms of a multidimensional torus fail to commute
and so typically a finite set of automorphisms generate a rather large
group. An extreme case is the study of actions of groups of
automorphisms which are commensurable to an arithmetic lattice in a
higher rank Lie group. Certainly in this case one expects very few
closed invariant sets and invariant measures, but much of the rigidity
comes from the group itself and not from the particular action at hand.
The study of actions of these ``big'' groups is a deep and fascinating topic in its own
right, but has a completely different flavor and is beyond the scope of these lectures.

Here we focus on the other extreme, the study of measures and closed
sets in $\T^d$ invariant under a multidimensional (i.e. not virtually
cyclic) Abelian group of toral automorphisms.
Certainly there is nothing rigid about the way $\Z^n$ acts per se:
virtually anything which is known about $\Z$ actions is also known for
$\Z^n$ actions: indeed for actions of amenable groups
\cite{Ornstein-Weiss-amenable-groups}. The surprising fact is that many natural $\Z^n$ or $\R^n$ actions have (at least conjecturally)
this unexpected scarcity of closed invariant sets and invariant
measures.

Some care is needed when stating the exact conditions for rigidity. For
example, consider the $\Z^2$-action $\alpha$ generated by the two
commuting toral automorphisms on $\T^4$ given by the integer matrices:
\begin{equation}
A = \begin {pmatrix}
2 & 1 & 0 & 0 \\
1 & 1 & 0 & 0 \\
0 & 0 & 1 & 0 \\
0 & 0 & 0 & 1
\end {pmatrix} , \qquad
B =
\begin {pmatrix}
1 & 0 & 0 & 0 \\
0 & 1 & 0 & 0 \\
0 & 0 & 2 & 1 \\
0 & 0 & 1 & 1.
\end {pmatrix}
\end{equation}
This action is a faithful $\Z^2$-action (i.e. no $\mathbf n \in \Z^2$
acts trivially) but it is clear no measure or topological rigidity can
hold here since the projection $\pi : \T^4 \to \T^2$ obtained by
omitting the two last coordinates is a factor map  that intertwines the
$\Z^2$-action $\alpha$ on $\T^4$ with a $\Z$ action on $\T^2$. There is
no shortage of closed invariant sets and invariant measures for any
$\Z$ action by a hyperbolic toral automorphism on $\T^2$, and any such
invariant set or invariant measure can be lifted to an invariant set or
measure respectively on $( \T^4 , \alpha )$. It is a conjecture of
Katok and Spatzier that this type of construction is the only source of
non-algebraic $\alpha$ ergodic measures for a $\Z^d$ action on a torus
with $d > 1$.

A $\Z^d$ action $\alpha$ on $\T^n$ is said to be {\bf irreducible} if
the only closed, $\alpha$-invariant proper subgroups of $\T^n$ are
finite. Such actions do not have any lower dimensional factors (since
otherwise the kernel of the factor map will be a closed infinite proper
subgroup). It is {\bf strongly irreducible} if the same holds for every
finite index subgroup of $\Z^d$. For simplicity we assume that $\alpha$
is a faithful, strongly irreducible $\Z^d$ action on $\T^n$ for some $d
> 1$, and then the above method of creating strange invariant sets and
measures is not applicable, though even if one understands well the
irreducible case, the passage to the general case involves new issues
(see e.g. \cite{Kalinin-Katok, Katok-Spatzier,
Katok-Spatzier-corrections, Meiri-multidimensional, Meiri-Peres}). We
also mention that some of the most interesting applications have
involved considering non irreducible actions \cite{Kalinin-Katok,
Katok-Katok-Schmidt}.

Using the same basic strategy as Furstenberg, D. Berend extended
Furstenberg theorem to the context of  semigroups of commuting
endomorphisms of a torus and more general abelian groups
\cite{Berend-minimal-tori, Berend-invariant-groups,
Berend-invariant-tori}.

As mentioned earlier, Katok and Spatzier generalized Rudolph's proof to
this context (as well as to the type of examples considered in the next
section). As before, one needs an assumption about positive entropy of
the $\Z$ action generated by $\alpha_{ \mathbf n }$ for some $\mathbf n
\in \Z^d$. But in most cases, an ergodicity assumption is also needed.
We only discuss the case of irreducible actions, though Katok and
Spatzier treat also the general case.

\begin{Theorem}[Katok and Spatzier \cite{Katok-Spatzier,
Katok-Spatzier-corrections}]
Suppose $\alpha$ is a strongly irreducible $\Z^d$ faithful action on
$\T^n$ with $d > 1$. Let $\mu$ be an $\alpha$-invariant and ergodic
measure on $\T^n$, and suppose that for some $\mathbf n \in \Z^d$, the
entropy of $\mu$ with respect to the $\Z$ action generated by $\alpha_{
\mathbf n }$ is positive. Suppose further that one of the following
holds:
\begin{enumerate}
\item the $\Z^d$ action $\alpha$ on $\T^n$ is totally non symplectic
(we do not provide the exact definition here, but for example a
faithful $\Z^d$ action on $\T^{ d + 1 }$ is always of this kind);
\item the action given by $\alpha$ on the measure space $( X , \Borel ,
\mu )$ satisfies some ergodicity conditions; a natural condition to use

is that $\alpha$ is weakly mixing as a $\Z^d$ action which in
particular implies that for any $\mathbf n \in \Z^d$, we have that
$\mu$ is ergodic for the $\Z$ action generated by $\alpha_{ \mathbf n
}$.
\end{enumerate}
Then $\mu$ is Lebesgue measure.
\end{Theorem}

Jointly with M. Einsiedler, we have been able to
eliminate the need for assuming any ergodicity condition or total non
symplecticity.

\begin{Theorem} [Einsiedler and L. \cite{Einsiedler-Lindenstrauss-ERA}] \label{Einsiedler-Lindenstrauss theorem}
Suppose $\alpha$ is a strongly irreducible $\Z^d$ faithful action on
$\T^n$ with $d > 1$. Let $\mu$ be an $\alpha$-invariant and ergodic
measure on $\T^n$, and suppose that for some $\mathbf n \in \Z^d$, the
entropy of $\mu$ with respect to the $\Z$ action generated by $\alpha_{
\mathbf n }$ is positive. Then $\mu$ is Lebesgue measure.
\end{Theorem}

We remark that the techniques of \cite{Einsiedler-Lindenstrauss-ERA} apply equally well to the reducible case, but necessarily the statement in this case is more cumbersome. We refer to \cite{Einsiedler-Lindenstrauss-ERA} for additional details.

The proof of Theorem~\ref{Einsiedler-Lindenstrauss theorem} uses a new observation regarding measures invariant under a {\bf single}, possibly hyperbolic, toral automorphism. In order to present this observation, we first needs to explain how one can attach for any sub foliation of the unstable foliation its contribution to the entropy.

Let $T:x \to Ax \bmod \Z ^ n$ be such an irreducible toral automorphism, with $A \in \GL (n, \Z)$,  and $\mu$ a $T$-invariant measure.
Let $\Phi ^ +$ be the set of eigenvalues of $A$ with absolute value $> 1$. To any nonempty $S \subset \Phi^+$ which is closed under complex conjugation we set $W_S \subset \R ^ n$ to be the space spanned by the corresponding eigenvectors, and the orbits $x+W_S$ of this space form a $T$-invariant expanding foliation of $\T ^ n$. In this case one can define an  entropy contribution  $h_{ \mu}(T, W_S)$  which roughly corresponds to the part of the total metric (i.e. measure theoretic) entropy $h_\mu (T)$ of $T$ which can be detected using the conditional measure $\mu$ induces on $W_S$-leaves, and these satisfy the following properties:
\begin{enumerate}
\item for any $S$ as above, $h_\mu (T) \geq h_\mu (T, W_S)$, with equality for $S = \Phi ^ +$.
\item $0 \leq h _ \mu (T, W _ S) \leq \sum_ {s \in S} \log \absolute s$; for Lebesgue measure $\lambda$, $h _ \lambda (T, W _ S) = \sum_ {s \in S} \log \absolute s$.
\end{enumerate}
As is often the case in ergodic theory, one is also interested in the corresponding relative quantities: if $\mathcal{A}$ is a $T$-invariant sigma algebra of Borel subsets of $\T ^ n$, one can define the conditional entropies $h_\mu (T | \mathcal{A})$ and $h_\mu (T, W_S | \mathcal{A})$ which have similar properties.
Precise definitions of all these quantities and additional references can be found in \cite [Sect. 3]{ Einsiedler-Lindenstrauss-ERA}.

The basic inequality at the heart of the approach of \cite{Einsiedler-Lindenstrauss-ERA} is the following:
\begin{Theorem} [ {\cite [Thm. 4.1]{ Einsiedler-Lindenstrauss-ERA}}] \label{theorem: entropy inequality}
Let $T$ be an irreducible toral automorphism, and $\mu$ a $T$-invariant measure. Let $W_S \subset \R ^ n$ be an expanding $T$-invariant subspace as above. Then
\begin{equation} \label{entropy inequality}
h _ \mu (T, W _ S) \leq \frac {h _ \lambda (T, W _ S) }{  h _ \lambda (T)} h _ \mu (T).
\end{equation}
More generally, for every $T$-invariant sigma algebra $\mathcal{A}$ it holds that $h _ \mu (T, W _ S \cond \mathcal{A}) \leq \frac {h _ \lambda (T, W _ S) }{  h _ \lambda (T)} h _ \mu (T \cond \mathcal{A})$
\end{Theorem}
We emphasize that \equ{entropy inequality} is false for reducible $T$. It is an interesting question when can equality be attained.

Returning to the case of an irreducible $\Z ^ d$ action $\alpha$ for $d > 1$, in conjunction with a lemma from \cite{Einsiedler-Katok}, Theorem~\ref{theorem: entropy inequality} implies that for any $ \alpha$-invariant $\mu$ and $\mathcal{A}$ there is a single constant $c_{\mu, \mathcal{A}} \in [ 0, 1 ]$ such that for any $\mathbf n \in \Z ^ d$
\begin{equation} \label{equation: one dimensionality of entropy}
h _ \mu (\alpha _ \mathbf n \cond \mathcal{A}) = c _ {\mu, \mathcal{A}} h _ \lambda (\alpha _ \mathbf n)
.\end{equation}
This one dimensionality of the entropy function for $\alpha$-invariant measures can be expected in the one-dimensional case covered by Rudolph's theorem (Theorem~\ref{Rudolph theorem}), where entropy is closely related to Hausdorff dimension of $\mu$, and indeed \equ{equation: one dimensionality of entropy} plays a prominent role in Rudolph's proof (see \cite [Thm. 3.7]{ Rudolph-2-and-3}). That it holds also in the higher dimensional irreducible case is quite surprising, and allows one to prove Theorem~\ref{Einsiedler-Lindenstrauss theorem} along the lines of Rudolph's original proof.\footnote{The alert reader may notice that in Theorem~\ref{Einsiedler-Lindenstrauss theorem} strong irreducibility is assumed while in Theorem~\ref{theorem: entropy inequality} no such assumption is made. The reason is that in Rudolph's proof strong irreducibility is implicitly assumed; indeed Theorem~\ref{Einsiedler-Lindenstrauss theorem} as stated is false in the irreducible but not strongly irreducible case.}

Finally, we mention that there has been some work lately on classifying
closed subsets of $\T^n$ invariant under an action by more general
semigroups of automorphisms (or endomorphisms) of the torus
\cite{Muchnik-Zariski-dense, Starkov-toral-automorphisms,
Starkov-toral-automorphisms-II, Muchnik-toral-actions}. All of these
results assume in particular that the action is strongly irreducible
and the group acting is not virtually cyclic; under these assumption
\cite{Muchnik-toral-actions} provides a complete classification.

\section{rigidity on homogeneous spaces}

Here we present another class of highly interesting examples where
topological and measure rigidity of multiparameter actions is expected
to hold.

One of the simplest, and arguably the most interesting cases is the
following: let $X = \SL ( 3 , \Z ) \backslash \SL ( 3 , \R )$. This is
simply the space of lattices in $\R^3$ of covolume one. It is not a
compact space: if one takes a sequence of lattices $\Lambda_i$ in
$\R^3$ with each $\Lambda_i$ containing a vector $v_i$ with $\norm {
v_i } \to 0$ then no subsequence of $\Lambda_i$ has a limit in $X$.
However, $X$ has a (unique) probability measure which is invariant

under the natural action of $\SL ( 3 , \R )$ on $X$. The diagonal
matrices
\begin{equation*}
H = \left\{ \begin {pmatrix}
e^t & 0 & 0 \\
0 & e^s & 0 \\
0 & 0 & e^{ - s - t }
\end {pmatrix}
: s , t \in \R \right\} \subset \SL ( 3 , \R )
\end{equation*}
give a very nice $\R^2$ action on $X$ by multiplication from the right.

We mention in the passing that for any one parameter subgroup of $H$
there is a big collections of closed invariant sets and of invariant
measures; but these collections are much less understood than in the
case of (hyperbolic) toral automorphisms. We also mention that there
are known irregular $H$-invariant closed subsets of $X$ arising, as in
the previous section, from special situations in which the action is
essentially reduced to an action of a one parameter group. For this
particular $X$ this can happen if we are considering the orbit closure
of a point $\Lambda \in X$ under the $H$-action when there is
one one parameter subgroup of $H$, say $\left\{ a ( t ) : t \in \R
\right\}$ so that $a ( t ) \Lambda \to \infty$ as $\absolute
t \to \infty$. For more general examples of this kind, such reduction
to an action of a one parameter group can happen even in the compact
case. Note that when classifying invariant probability measures, by
Poincare recurrence, the particular situation described above cannot
happen, though in general, the same kind of difficulties exist when
classifying invariant probability measures (there is a instructive
example of this phenomenon due to M. Rees \cite{Rees-example}).

To avoid this complication of possible bad behavior of certain
diverging orbits, we can consider the following partial classification
of $H$ orbit on $X = \SL ( 3 , \Z ) \backslash \SL ( 3 , \R )$. We note
that even for this specific $X$, this is only a special case of a more
general conjecture regarding closed $H$-orbits.

\begin{Conjecture}[G. Margulis \cite{Margulis-Oppenheim-conjecture}]
\label{Margulis conjecture in special case}
Any compact $H$-invariant subset of $\SL ( 3 , \Z ) \backslash \SL ( 3
, \R )$ is a union of compact $H$-orbits.
\end{Conjecture}

It turns out that even this very special case of topological rigidity
for homogeneous spaces implies the following long-standing conjecture
of Littlewood:

\begin{Conjecture}[Littlewood (c. 1930)]
Let $\norm x$ denote the distance from $x \in \R$ to the closest
integer. Then
\begin{equation} \label{equation regarding Littlewood}
\liminf_{ n \to \infty } n \norm { n \alpha } \norm { n \beta } = 0
\end{equation}
for any real numbers $\alpha$ and $\beta$.
\end{Conjecture}

This implication has been discovered in a different terminology long
before Furstenberg's pioneering work regarding the rigidity of
multiparameter actions by Cassels and Swinnerton-Dyer
\cite{Cassels-Swinnerton-Dyer}; however, it was Margulis who first
recast this in dynamical terms \cite{Margulis-Oppenheim-conjecture}.
We note that the current state-of-the-art regarding Conjecture~\ref{Margulis conjecture in special case}, Theorem~\ref{Einsiedler-Katok-Lindenstrauss theorem} below, while not sufficient to imply Littlewood's conjecture, does imply the following nontrivial estimate towards this conjecture:

\begin{Theorem} [ Einsiedler, Katok and L. \cite{Einsiedler-Katok-Lindenstrauss} ] \label{theorem regarding Littlewood}
The set of $(\alpha, \beta) \in \R ^2$ for which \equ{equation regarding Littlewood} is not satisfied has Hausdorff dimension zero.
\end{Theorem}
A weaker result, namely that this Hausdorff dimension is at most one, can be derived from Theorem~\ref{theorem of Einsiedler and Katok}.

The most explicit form of a general rigidity conjecture for $\R^d$ and
more general actions on homogeneous spaces has been given by Margulis
\cite{Margulis-conjectures}. A less explicit conjecture in the same
spirit has been given by Katok and Spatzier in \cite{Katok-Spatzier}.

We take $G$ to be a connected Lie group, $\Gamma$ a lattice in $G$,
i.e. discrete subgroup of $G$ of finite covolume, and $H < G$ a
connected subgroup. One of the highlights of the theory of flows on
homogeneous spaces is the classification of closed invariant sets and
invariant measures on $\Gamma \backslash G$ when $H$ is generated by
unipotent one parameter subgroups
by M. Ratner \cite{Ratner-Annals, Ratner-Duke}; some years
earlier an important special case allowed Margulis to resolve the
long-standing Oppenheim conjecture (\cite{Margulis-Oppenheim-original}; see \cite{Dani-Margulis} for a very elegant and accessible treatment as well as the survey \cite{Margulis-Oppenheim-conjecture} for more details); we also mention
that this measure and topological rigidity of unipotent flows was
conjectured by Dani and Raghunathan respectively. A valuable recent addition to the literature on this deep and interesting subject is  \cite{Witte-Ratner-theorem}.\footnote{The simplest case is the horocycle flow on compact surfaces of constant negative curvature, for which minimality (i.e. nonexistentance of closed invariant subsets) is classical and due to Hedlund \cite{Hedlund-dynamics}; unique ergodicity of this flow, which was the first result regarding measure rigidity of unipotent flows, has been proved by Furstenberg in \cite{Furstenberg-horocycle}.}
\nocite{Margulis-Tomanov}

Taking $H$ to be as in Conjecture \ref{Margulis conjecture in special
case} is perhaps a prototypical example of the kind of groups not
covered by Ratner's theorem. In his conjecture, Margulis
only assumes that $H$ is generated by elements $h \in H$ so that all
eigenvalues of $Ad_G ( h )$ are real. Margulis gives a conjecture both
in the topological category and in the measure theoretic category:

\begin{Conjecture}[Margulis \cite{Margulis-conjectures}]
\label{Margulis conjecture about orbits}
Let $G , \Gamma$ and $H$ be as above. For any $x \in \Gamma \backslash
G$ one of the following holds:
\begin{enumerate}
\item $\overline { Hx }$ is homogeneous, i.e. is a closed orbit of a
closed connected subgroup $L < G$
\item there exists a closed connected subgroup $F < G$ and a continuous
epimorphism $\phi$ of $F$ onto a Lie group $L$ such that $H < F$, $Fx$
is closed in $\Gamma \backslash G$, $\phi ( F_x )$ is closed in $L$
where $F_x$ denote the stabilizer $\left\{ g \in F : gx = x \right\}$,
and $\phi ( H )$ is a one parameter subgroup of $L$ containing no
nontrivial $Ad_L$-unipotent elements.
\end{enumerate}
\end{Conjecture}

\begin{Conjecture}[Margulis \cite{Margulis-conjectures}]
\label{Margulis conjecture about measures}
Let $G , \Gamma$ and $H$ be as above. Let $\mu$ be an $H$-invariant
$H$-ergodic measure on $\Gamma \backslash G$. Suppose that for every $x
\in \supp \mu$ the statement (2) in the formulation of Conjecture
\ref{Margulis conjecture about orbits} does not hold. Then $\mu$ is
algebraic (i.e. the invariant measure on the closed orbit of some
subgroup $L < G$ containing $H$).
\end{Conjecture}

It seems to me that the second conjecture is overly restrictive. For
the record, we state the following modification of Conjecture
\ref{Margulis conjecture about measures}:

\begin{Conjecture}
Let $G , \Gamma$ and $H$ be as above. Let $\mu$ be an $H$-invariant
$H$-ergodic measure on $\Gamma \backslash G$. Suppose that for
$\mu$-almost every $x$ the statement (2) in the formulation of
Conjecture \ref{Margulis conjecture about orbits} does not hold. Then
$\mu$ is algebraic.
\end{Conjecture}

Note that these conjectures give no information for $H$ a one-parameter
diagonalizble group; and despite the fact that they treat non
commutative $H$ the most interesting case seems to be $H$ commutative
(see discussion and remarks following Conjecture 2 in
\cite{Margulis-conjectures}). From now on we will only consider actions
of commutative $H$.

Unlike the case for the original $\times _ n, \times _ m$-problem and its extensions discussed in the previous section, the study of invariant measures seems to have fared better than that of the topological questions.
Regarding the topological question, S. Mozes \cite{Mozes-quaternions} has proved the following: let $G=\PGL (2, \Q _ p) \times \PGL (2, \Q _ q)$, $H<G$ the product of the groups of diagonal matrices in the two factors of $G$, and  $\Gamma$ an irreducible lattice. Then if  $\overline { Hx }$ contains a compact
$H$ orbit,
then $\overline { Hx }$ is homogeneous. This result has a nice interpretation in terms of certain two-dimensional tilings. Jointly with Barak Weiss \cite{Lindenstrauss-Barak} we have proved a similar result for $H$ the full diagonal group in $G = \SL ( n , \R )$ and $\Gamma$ e.g. $ \SL ( n , \Z )$; a substantial complication is the existence of intermediate groups $H<L<G$ with closed orbits for $n$ composite. Generalization to general $G$ seems to require additional ideas.
We also mention a recent paper
by G. Tomanov and Barak Weiss \cite{Tomanov-Weiss} classifying closed
$H$-orbits in $\Gamma \backslash G$ for algebraic groups and $H$ a
maximal $\R$-split torus (e.g. for $G = \SL ( n , \R )$, $H$ is the
group of positive diagonal matrices). Their treatment generalizes an
unpublished result of Margulis (see the appendix to
\cite{Tomanov-Weiss}).
Finally, we remark that the latest results towards the measure theoretic conjecture have substantial implications regarding orbit closures.

We now turn to the measure theoretic conjecture.
In their papers, Katok and Spatzier treat simultaneously the
homogeneous case and the case of multidimensional action of toral
automorphisms. However (at least when $G$ is semi simple and $H$ is a
multidimensional subgroup of the Cartan group of $G$) this action never
satisfies their total non symplecticity assumption. Because of this, in
order to deduce that an $H$ invariant measure is homogeneous both an
ergodicity assumption and an assumption regarding entropy is needed.

The entropy assumption is a major disadvantage. While the Katok-Spatzier method and its extensions has been used very effectively to prove isomorphism rigidity in many cases (i.e. that any measurable isomorphism between two multiparameter actions is algebraic; for a study of isomorphism rigidity in the locally homogeneous context we refer the reader to the recent \cite{Kalinin-Spatzier} by B. Kalinin and Spatzier), its applicability to number theoretic questions is extremely limited.
The reason for this is that typically, when applying measure rigidity to prove e.g. equidistribution, one considers the weak$^{*}$ limit of a sequence of measures $ \mu _ i$.
While entropy is usually well behaved under weak$^{*}$ limits, mixing properties such as ergodicity are not.
This principle is well illustrated in e.g. the proof of Theorem~\ref{theorem regarding Littlewood} in \cite [Part 2]{ Einsiedler-Katok-Lindenstrauss} (see also \cite{Johnson-Rudolph} for a related discussion) and is implicit in \cite{Bourgain-Lindenstrauss, Lindenstrauss-quantum}.

Recently, two rather different and complimentary methods to overcome the need for an ergodicity assumption for the local homogeneous context have been discovered (yet another idea, which to date has not been used in this context but is very useful for e.g. commuting automorphisms of tori has been briefly exposed in the previous section).
Chronologically the first of these has been by
Einsiedler and Katok \cite{Einsiedler-Katok}; motivated by this result and by a beautiful and patient exposition of Ratner's work on unipotent flows by Dave Witte Morris (who has later also written a book on this topic \cite{Witte-Ratner-theorem}) I was fortunate enough to find another completely different approach shortly thereafter \cite{Lindenstrauss-quantum}.

In their paper, Einsiedler and Katok deal with general $\R$-split group $G$ with rank $> 1$; for simplicity (and also because the result is
strongest in this case) we specialize to the case of $G = \SL ( n , \R )$.
In this case Einsiedler and Katok prove the following:

\begin{Theorem}[Einsiedler and Katok \cite{Einsiedler-Katok}]
\label{theorem of Einsiedler and Katok}
Let $G = \SL ( n , \R )$ and let $\Gamma<G$ be discrete. Let $H <
G$ be the subgroup of positive diagonal matrices. Let $\mu$ be a
$H$-invariant and ergodic measure on $\Gamma \backslash G$. Assume that
the entropy of $\mu$ with respect to all one parameter subgroup of $H$
is positive. Then $\mu$ is the $G$ invariant measure on $\Gamma
\backslash G$.
\end{Theorem}

The proof of this result uses heavily the non commutativity of the one
parameter unipotent subgroups of $G$ normalized by $H$; such use is
already hinted at in \cite{Katok-Spatzier}. But it also employs a new
observation regarding how a measure invariant under a multiparameter
action decomposes when restricted to subfoliations of the stable
foliation. This ingredient has already found several additional
applications, including in \cite{Einsiedler-Lindenstrauss-ERA}.

Using a completely different method, closely related to Ratner's
techniques for proving the Raghunathan conjectures, the following (and
more general statements) can be proved (see the discussion following the closely related Theorem~\ref{theorem for que}). We note that the
Einsiedler-Katok method gives very little in this case, and that the
entropy assumption is weaker than that of Theorem \ref{theorem of
Einsiedler and Katok}.

\begin{Theorem} [L. \cite{Lindenstrauss-quantum}]
Let $G_1 = G_2 = \SL ( 2 , \R )$, $G = G_1 \times G_2$, and $\Gamma$ an
irreducible lattice in $G$ (i.e. a lattice in $G$ whose projections to
$G_1$ and $G_2$ are dense). Let $H$ be the product of the group of
positive diagonal matrices in $G_1$ with the group of positive diagonal
matrices in $G_2$. Suppose $\mu$ is $H$ invariant and ergodic, and
suppose that the entropy of $\mu$ with respect to the action of the one
parameter group of positive diagonal matrices in $G_1$ is positive.
Then $\mu$ is the $G$ invariant measure on $\Gamma \backslash G$.
\end{Theorem}

It turns out that the argument of \cite{Lindenstrauss-quantum} is applicable for $\Gamma \backslash \SL (n, \R)$ precisely for positive entropy measures for which the argument of \cite{Einsiedler-Katok} does not apply(!); combining these arguments one proves the following strengthening\footnote{the restriction to $\SL (n, \Z)$ is not a technicality: there are lattices, even cocompact latices, which have non algebraic $H$ invariant measures with positive entropy with respect to a one parameter subgroup of $H$ (indeed, the Rees examples alluded to earlier can have positive entropy). Note also that even for $\SL (n, \Z)$, under the weaker assumptions of Theorem~\ref{Einsiedler-Katok-Lindenstrauss theorem}, $\mu$ need not be the $G$ invariant measure on $\Gamma \backslash G$, but can be an algebraic measure on a proper subspace.} of Theorem~\ref{theorem of Einsiedler and Katok}:

\begin{Theorem} [ Einsiedler, Katok and L. \cite{Einsiedler-Katok-Lindenstrauss}] \label{Einsiedler-Katok-Lindenstrauss theorem}
Let $G = \SL ( n , \R )$ and $H <
G$ be the subgroup of positive diagonal matrices. Let $\mu$ be a
$H$-invariant and ergodic measure on $\SL (n, \Z) \backslash G$. Assume that there is some one parameter subgroup of $H$ with respect to which
$\mu$ has positive measure. Then $\mu$ is algebraic, and is not compactly supported. If $n$ is prime then $\mu$ is the $G$ invariant measure on $\SL (n, \Z) \backslash G$.
\end{Theorem}

\section{Recurrent measures}

We now return for a moment to $\T$. B. Host gave a proof of Rudolph's theorem regarding invariant measures on $\T$ along the following
outline: first he showed that every $\times_n$ invariant measure $\mu$
all of whose $\times_n$-ergodic components have positive entropy
are in an appropriate sense (which will be defined momentarily) recurrent under the
action of the {\bf additive} group
\begin{equation*}
N_n = \left\{ \frac{ q }{ n^k } : k \N , q \in \Z / ( n^k ) \Z
\right\}.
\end{equation*}
Host then shows that while there is an abundance of measures on $\T$
invariant under $\times_m$, if $m$ and $n$ are relatively prime the
only $\times_m$ invariant measure which is $N_n$-recurrent is
Lebesgue measure.

In fact, Host proved more: if one assumes that $\mu$ is
$\times_n$-invariant, and the entropy of the $\times_n$-ergodic
components of $\mu$ are bounded from below by some positive number then
$\mu$ satisfies a stronger quantitative condition that implies
$N_n$-recurrence, and that for any measure $\mu$ which satisfies this condition, the orbit of $\mu$-almost every $x$ under successive
iterates of $\times_m$ is uniformly distributed with respect to the
Lebesgue measure on $\T$ (which by the pointwise ergodic theorem
clearly implies that if $\mu$ is also $\times_m$ invariant then
$\mu$ is Lebesgue measure). Extensions of Host's result to the
multidimensional case can be found in \cite{Host-multidimensional,
Meiri-multidimensional};  other related papers that in particular treat
the case of $m , n$ which are not relatively prime can be found in
\cite{Lindenstrauss-Meiri-Peres, Lindenstrauss-p-adic}.

The simplest case is when considering an action of a locally compact
group $M$ on a complete separable metric space $X$. We say that the
measure $\mu$ on $X$ is $M$-recurrent if for every subset $B \subset
X$ with $\mu ( B ) > 0$
we have that for almost every $x \in B$, for every compact $K \subset M$,
there is a $g \in M \setminus K$ so that $gx \in B$, in other words,
if Poincare recurrence holds for this
action.  This certainly seems the very least one needs to assume about a measurable action to get nontrivial dynamics.

In the very special case of $\Z$ or $\R$ actions and $\mu$ quasi-invariant (i.e. while $\mu$ may not be conserved by the action, the measure class is), what we call recurrence is known as conservativity and plays an important
role; for example, see \cite[Sect. 1.1]{ Aaronson-book}. Host may well be the first to have made essential use of this condition when even the measure class is not invariant; he is certainly the first to have used it in  the context surveyed in this paper.

It is customary in the literature to refer to a space $X$ on which a
group $G$ acts as a $G$-space (to avoid pathologies, it is helpful to
assume e.g. that $G$ acts freely on a dense subset of $X$).  A
generalization of this is the notion of a
$(G,T)$-space, for a nice\footnote{locally compact separable and metrizable} group $G$, and $T$ a homogeneous space for $G$ (i.e. a locally compact metric space $T$ on which $G$ acts transitively)  with a
distinguished point $e \in T$.

\begin{Definition}

\label{definition of t-space}
A locally compact metric space $X$ is said to be a
$(G, T)$-space if there is some open cover
$\mathfrak T $ of $X$ by relatively compact sets, and for
every $U \in \mathfrak T$ and $x \in U$ a map $t _ {U, x}: T \to X$ with the following
properties:
\renewcommand{\theenumi}{A-\arabic{enumi}}
\begin{enumerate}
\item For every $x \in U \in \mathfrak T$, we have that $t _
{U,x}(e) =
x$.
\item The map $ (x, t) \mapsto t _ {U, x} (t)$ is a continuous map $U \times T \to X$.
\item \label{item describing equivalence classes} For any $x \in
U \in \mathfrak T$, for any $y$ in the image $\image t _ {U, x}$ of $t _ {U, x}$
and $V \in \mathfrak T$ containing $y$, there is
a $\theta \in G$  so that
\begin{equation}
\label{equation for Greek theta}
t _ {V,y} \circ \theta =
t _ {U,x}.
\end{equation}
In particular, for any $x,y,U,V$as above $\image t _ {U, x} = \image t _ {V, y}$ and $ t _ {U, x}$ is injective if and only if
$t _ {V, x}$ is.
\item for any $U \in \mathfrak T$ there is some open neighborhood $B ^ U$ of the distinguished point $e \in T$  so that $t _ {U, x}$ is injective on $\overline {B ^ U}$ for all $x \in
U$.
\item For a dense set of $x \in X$,  for some (equivalently for
every) $U \in \mathfrak T$ containing $X$ the map $t _ {U, x} $
is injective.
\end{enumerate}
\renewcommand{\theenumi}{\arabic{enumi}}
\end{Definition}

We will be only interested in the case where $G$ acts by isometries, in other words when $G$ is a subgroup of the group $\isomof T$ of isometries of $T$. For simplicity, we shall call a $(\isomof T, T)$-space simply a $T$-space.

In less technical terms, a $T$-space is a space which is foliated in a
continuous way by leaves, such that on every leaf there is a well
defined metric, and every leaf is locally isometric to $T$ with respect
to this metric. Furthermore, for most leaves this is actually an
isometry between the leaf and $T$. A $(G,T)$-space for $G$ a proper subgroup of $\isomof T$ is a special kind of $T$-space with some additional structure.

\begin{Definition}
A measure $\mu$ on any $T$-space $X$ is said to be {\bf
$T$-recurrent} if for every $B \subset X$ with $\mu ( B ) > 0$, for
$\mu$-almost every $x \in B$ (and some $U \in \mathfrak{T}$ containing
$x$) for every compact $K \subset T$,
\begin{equation*}
t _ {U, x} (T \setminus K) \cap B \neq \emptyset
.\end{equation*}
\end{Definition}

If $\mu$ is a Radon\footnote{i.e. a regular Borel measure which is finite on compact sets} measure on a $(G,T)$-space, it induces on every $T$-leaf a
$\sigma$-finite measure which is well defined up to a multiplicative
scalar. We can either think of these measures as $\sigma$-finite but no longer Radon measures on $X$, or use the isomorphism between typical $T$-leaves and $T$ and consider these measures as nice Radon measures on $T$. We adopt here the latter point of view, and so formally these conditional measures are a system $\left\{ \mu ^ U _ {T, x} \right\} _ {U \in \mathfrak T, x \in U}$ of Radon measures on $T$ so that for every $U$, the map $ x \mapsto \mu ^{U}_{T,x}$ is Borel measurable\footnote{with respect to the weak$^{*}$ topology on the space of Radon measures on $T$ as the dual to the space of compactly supported continuous functions} and so that there is a set of full measure $X' \subset X$ so that for every $x, y \in X '$ and $U,V \in \mathfrak T,g \in G$ such that $t _ {U, x} \circ g = t _ {V, y}$ then
\begin{equation*}
\mu ^ U _ {T, x} = g_*\mu ^ V _ {T, y}
\end{equation*}
i.e. so that formally $(t _ {U, x}) _ {*} \mu ^ U _ {T, x}$ is equal to $(t _ {V, y}) _ {*} \mu ^ V _ {T, y}$.
Defining these measures carefully is slightly tricky (see \cite [Sect. 3]{ Lindenstrauss-quantum}; \cite{Kalinin-Katok-Seattle, Lindenstrauss-Chile} may also be useful).

These conditional measures reflect the dynamics of $\mu$ with respect to the foliation: for example, if $T$ is a group and the $T$-structure comes from some action of $T$, then $\mu$ is $T$-invariant if and only if $\mu _ {T,x}$ is Haar measure for almost every $x$ (the dependence on the set $U$ in the atlas $\mathfrak T$ is trivial and therefore omitted in this case). Furthermore, and this is an important fact, if $\mu$ is a probability measure on a $(G, T)$-space then $\mu$ is
$T$-recurrent if and only if $\mu ^ U _ {T, x}$ is an infinite measure a.s.\footnote{for $\Z$-action which preserve
the measure class of $\mu$, this result is known as the Halmos Recurrence Theorem (see \cite[Sect. 1.1]{ Aaronson-book})}

Jointly with Klaus Schmidt, and using a method reminiscent (but distinct) of
that of Host, we have proved

\begin{Theorem}[Schmidt and L. \cite{Lindenstrauss-Schmidt}]
Let $\alpha$ be a nonexpansive, ergodic and totally irreducible automorphism of a compact connected abelian group $X$  (for example, let $\alpha$ be a nonhyperbolic toral automorphism). And let $\mu $ be an $\alpha $-invariant probability measure on $X$. Then $\mu$ is recurrent with respect to the central foliation of $X$ if and only if it is Lebesgue measure on $X$.
\end{Theorem}

This theorem also implies that if, for example, $\mu$ is weakly mixing
for $\alpha$ then there is a Borel subset of $X$, which intersects
every central leaf in at most one point, and which has full
$\mu$-measure. A similar, but slightly weaker statement holds for a
general invariant measure which is singular with respect to Lebesgue
measure. Thus, in some sense the only invariant measure on $X$ which is
truly nonhyperbolic is Lebesgue measure.

We contrast this with the situation of another partially hyperbolic
dynamical system: the time one map $g_1$ on a negatively curved compact
manifold, where {\bf all} $g_1$-invariant measures are recurrent
with respect to the central foliation, which in this case is simply the
foliation of the unit tangent bundle of the manifold into orbits of the
geodesic flow.

\section{Rigidity of Hecke recurrent measures on arithmetic surfaces}
\label{section: rigidity on arithmetic surfaces}

A similar theorem can be proved in the context of homogeneous spaces,
and can be used to settle an important special case of a conjecture of
Rudnick and Sarnak \cite{Rudnick-Sarnak} regarding the behavior of
eigenfunctions of the Laplacian on negatively curved manifolds.

Take $X = \Gamma \backslash \PGL ( 2 , \R )$, where $\Gamma$ is either
a congruence subgroup of $\PGL ( 2 , \Z )$ or of certain lattices that
arise from quaternionic division algebras over $\Q$ that are unramified
over $\R$. The later lattices are slightly harder to define but have
the advantage that $X$ is compact.

In both cases, for all but finitely many prime $p$, there is a map
$T_p$ from $X$ to $( p + 1 )$-tuples of points of $X$ called the {\bf
Hecke correspondence}. For example, if $\Gamma = \SL ( 2 , \Z )$, one
can define $T_p$ of a point $x \in X$ as follows: choose some $g \in
\SL ( 2 , \R )$ so that $x = \Gamma g$. Then
\begin{equation*}
T_p ( x ) = \left\{ \Gamma \twobytwo p 0 0 1 g , \Gamma \twobytwo p 0 1
1 g , \dots , \Gamma \twobytwo p 0 { p - 1 } 1 g , \Gamma \twobytwo 1 0
0 p g \right\}.
\end{equation*}
Notice that while each individual element of the right hand side depend
on the choice of $g$, the set of $p + 1$ points depends only on $x$.
The above formula also shows that for any $x \in X , g \in \PGL ( 2 ,
\R )$
\begin{equation*}
T_p ( xg ) = T_p ( x ) g.
\end{equation*}

Using the Hecke correspondence we can define operators called the {\bf
Hecke operators} on $L^2 ( X )$ by
\begin{equation}
\label{equation with invariance of Hecke}
T_p ( f ) [ x ] = \sum_{ y \in T_p ( x ) } f ( y )
\end{equation}
which play a very important role in the spectral theory of $X$.
The Hecke correspondences have some additional nice properties, which
are not entirely obvious from the way we defined them here. For
example, if $y \in T_p ( x )$ then $x \in T_p ( y )$, which means that
the Hecke operators are self adjoint. Furthermore, for any $p$ and $q$,
the operators $T_p$ and $T_q$ commute.

Fix a prime $p$. We define the following equivalence relation on $X$:
$x \sim y$ if there is a chain of points $x_i$ with $x_0 = x$, \ $x_n =
y$ and $x_{ i + 1 } \in T_p ( x_i )$. We can give each equivalence
class a structure of a $p + 1$-regular graph if we connect every $x$
with every point in $T_p ( x )$, and this graph turns out to be a $p +
1$-regular tree for all $x$. By \equ{equation with invariance of Hecke}
this equivalence relation, and indeed even the structure of every
equivalence class as a $p + 1$-regular tree, is invariant under left
translation by elements of $\PGL ( 2 , \R )$. It is not hard to see
that if $T$ is a $p + 1$-regular tree, this partition of $X$ into $p +
1$-regular tree is a $T$-space.

Consider now the action of the one parameter group
$A=\left\{ a ( t ) = \twobytwo { e^t } 0 0 { e^{ - t } } \right\}$ on
$X$.
This action is equivalent to the geodesic flow on the hyperbolic
manifold $\Gamma \backslash \H$ with $\H$ the hyperbolic plane, which
is a prototypical example where measure rigidity does not hold.
However, it seems plausible that there will not be many measures $\mu$
which are invariant under $A$ and are $T$-recurrent. At the
current state of our knowledge, this still seems unapproachable. But if
one assumes in addition an entropy condition one can obtain the
following:

\begin{Theorem}[L. \cite{Lindenstrauss-quantum}]
\label{theorem for que}
Suppose $\mu$ is a $A$ invariant and $T$-recurrent measure on
$X$, where $T$ is a $p + 1$-regular tree as above. Suppose that
$\mu$-almost every $A$ ergodic component of $\mu$ has positive
entropy. Then $\mu$ is Lebesgue measure on $X$.
\end{Theorem}
We remark that as in Rudolph's theorem, a substantial source of
difficulty is that $\mu$ need not be ergodic under $A$.

Somewhat surprisingly, the main actor in the proof of Theorem~\ref{theorem for que} is not the action of the flow $A$ which preserves the measure, but rather of different flow: the horocycle flow given by the (right) action of the group
\begin{equation*}
U ^ + = \left\{ u (s) : s \in \R \right\} \qquad u (s) = \twobytwo 1 s 0 1
\end{equation*}
on $X$ (which a priori need not preserve neither $\mu$ nor even the measure class of $\mu$). Note that since $A$ normalizes the group $U^{+}$, the orbits of the group $U ^ +$ are preserved under $a (t)$ and are uniformly contracted by this flow, a relation which links these two actions in many subtle ways.

\begin{figure}
\begin{center}
\includegraphics[width=4.5in]{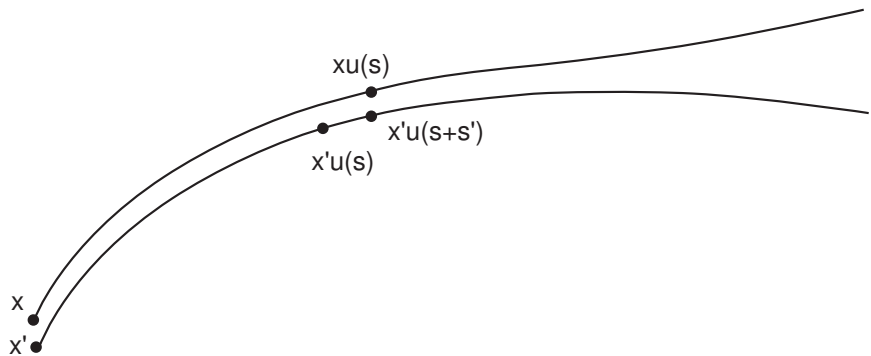}
\end{center}
\caption{Ratner's H-property}
\label{H figure}
\end{figure}

Indeed, the proof of Theorem~\ref{theorem for que} has many similarities to Ratner's earlier work on rigidity of unipotent flows \cite{Ratner-isomorphisms, Ratner-factors, Ratner-products}, and in particular uses Ratner's {\bf H-property} of horocycle flows (a related, but distinct property, the R-property, lies at the heart of Ratner's subsequent proof of the Raghunathan Conjecture; see \cite{Ratner-SL2} for more details). This property is a quantitative form of the following qualitative phenomenon: consider the trajectories under $U^{+}$ of two nearby points $x, x ' \in X$. Then there are two possibilities\footnote{for simplicity we assume $X$ compact; minor modifications are needed for the nonuniform case}:
\prefixitemswith{H-}
\begin{enumerate}
\item $xu (s), x ' u (s)$ remain close for all $s$, which can occur only if $x ' = x u (\eta)$ for some small $\eta \in \R$.
\item otherwise, there will be some smallest $s>0$ where $x u (s), x ' u (s)$ differ by a noticeable amount (say are of distance $\delta$ from each other, for some fixed but small $\delta$). Then most of the distance between $x u (s)$ and $x ' u (s)$ is along the direction of the $U ^ +$ flow, i.e. there is some $s '$ (roughly of order $\delta$) so that $x u (s)$ is close to $x ' u (s + s ')$ (see Figure \ref{H figure}).
\end{enumerate}
\prefixitemswith{}
The proof of the H-property is not difficult: it is merely an exercise in multiplying $2 \times 2$ matrices, but it has many far-reaching implications.

We now give a rough sketch of the proof of Theorem~\ref{theorem for que}, explaining how the H-property of the horocycle flow $U ^ +$ can be used to study measures invariant under the $A$-flow.

It follows from a slight generalization of a lemma of Einsiedler and Katok and a Fubini argument that if $\mu$ is $a(t)$-invariant, there is some set of full measure $X ' \subset X$ so that for every $x, y \in X '$ which are in the same Hecke tree (i.e. $x \sim y$) satisfy that $ \mu _ {x, U^+} = \mu _ {y,U^+}$. Using Hecke recurrence starting from almost any $x \in X '$ we can find a $ y \sim x$ in $ X '$ which is arbitrarily close to $x$.

By Luzin's theorem the map $z \mapsto \mu _ {z, U ^ +}$ is continuous in the weak$^*$ topology on a set of large measure, so since $y$ is nearby $x$ we can expect that  $ \mu _ {x, U^+}$ is fairly similar to $\mu _ {y, U^+}$. Of course this is no big surprise, and there are no contradiction so far since we know that in fact $\mu _ {x, U ^ +}$ is actually equal to $\mu _ {y, U ^ +}$. But now consider $yu (s)$ and $xu (s)$ which are still Hecke equivalent, for some $s$ large enough so that (H-2) above holds. Then by the properties of conditional measures typically $ \mu _ {xu (s), U ^ +}=\mu _ {yu (s), U ^ +}$. But now $xu(s)$ and $yu(s)$ are no longer close: indeed, $xu(s)$ is much closer to $yu(s+s')$ for some macroscopic $s'$. Thus we would expect $  \mu _ {xu (s), U ^ +}$ to be fairly similar to $\mu _ {yu (s +s'), U ^ +}= \mu _ {xu (s+s'), U ^ +}$. This hints we may be able to show that $\mu$ is $U^+$-invariant. Similarly one shows that $\mu$ is invariant under $U^-= \left\{ \twobytwo 1 0 s 1 \right\}$, so $\mu$ is Haar measure.

There are many tricky issues in making this crude idea work. One ingredient is a maximal inequality (sometimes known as a maximal ergodic theorem) for singular actions which appears in \cite[Appendix]{ Lindenstrauss-quantum} and is joint with D. Rudolph. This ergodic theorem was motivated by a similar idea in Rudolph's proof of Bernoullicity of the  Patterson-Sullivan measure \cite{Rudolph-Sullivan's-measure}.

\section {An application to quantum unique ergodicity}
\label{section: quantum unique ergodicity}

In this section we finally describe the quantum unique ergodicity conjecture of Rudnick and Sarnak, and explain how
it is related to Theorem \ref{theorem for que}.

\begin{Conjecture} [Rudnick and Sarnak \cite{Rudnick-Sarnak}]
\label{conjecture of Sarnak and Rudnick}
Let $M$ be a compact hyperbolic surface, and let $\phi_i$ a sequence of
linearly independent eigenfunctions of the Laplacian $\Delta$ on $M$,
normalized to have $L^2$-norm one. Then the probability measures
$\tilde \mu_i$ defined by $\tilde \mu_i ( A ) = \int_A \absolute {
\phi_i ( x ) }^2 d \haaron M ( x )$ tend in the weak$^*$ topology to
the uniform measure $\haaron m$.
\end{Conjecture}

A. I. {\v{S}}nirel$'$man, Y. Colin de Verdi{\`e}re and S. Zelditch have
shown in great generality (specifically, for any manifold on which the
geodesic flow is ergodic) that if one omits a subsequence of density 0
the remaining $\tilde \mu_i$ do indeed converge to $\haaron m$
\cite{Schnirelman, Colin-de-Verdiere-quantum,
Zelditch-uniform-distribution}. An important component of their proof
is the {\bf microlocal lift} of any weak$^*$ limit $\tilde \mu$ of a
subsequence of the $\tilde \mu_i$. The microlocal lift of $\tilde \mu$
is a measure $\mu$ on the unit tangent bundle $SM$ of $M$ whose
projection to $M$ is $\tilde \mu$, and most importantly it is always
invariant under the geodesic flow on $S M$. We shall call any measure
$\mu$ on $S M$ arising as a microlocal lift of a weak$^*$ limit of
$\tilde \mu_i$  a {\bf quantum limit}. Thus a slightly stronger form of
Conjecture \ref{conjecture of Sarnak and Rudnick} is the following
conjecture, also due to Rudnick and Sarnak:

\begin{Conjecture}[Quantum Unique Ergodicity Conjecture]
\label{quantum unique ergodicity conjecture}
For any compact hyperbolic surface $M$ the only quantum limit is the
uniform measure $\haaron { SM }$ on $S M$.
\end{Conjecture}

Consider now $M = \Gamma \backslash \H$, for $\Gamma$ one of the
arithmetic lattices considered above. If $\Gamma$ is a lattice
arising from a quaternionic division algebra over $\Q$, then $M$ is a
compact hyperbolic surface, precisely the kind of surface considered in
Conjectures \ref{conjecture of Sarnak and Rudnick} and \ref{quantum
unique ergodicity conjecture}.

If $\Gamma$ is a congruence sublattice of $\PGL ( 2 , \Z )$ then $M$
has finite volume, but is not compact.
It is not hard to modify conjectures \ref{conjecture of Sarnak and
Rudnick} and \ref{quantum unique ergodicity conjecture} so that they
cover also this case.
A special property of these very special surfaces of finite volume is
their spectral theory is not too far from that of compact surfaces of
the same area:
there is continuous spectrum (which of course does not exist for
compact surfaces), but it is very well understood, and as shown by
Selberg \cite{Selberg-trace-formula} the discrete spectrum of the Laplacian on these noncompact
surfaces obeys Weyl's law for compact surfaces of the same area. Since
Luo, Sarnak and Jakobson proved equidistribution for the continuous
spectrum \cite{Luo-Sarnak, Jakobson-equidistribution}, even for these
noncompact surfaces one only needs to study the weak$^*$ limits of
measures which arise from $L^2$-normalized eigenfunctions of the
Laplacian.

When looking at eigenfunctions of the Laplacian on the arithmetic
surfaces, it is natural to consider joint eigenfunctions both the
Laplacian and of all Hecke operators. Since the Hecke operators are
self adjoint operators that commute with each other and with the
Laplacian one can always find an orthonormal basis of the subspace of
$L^2 ( M )$ which corresponds to the discrete part of the spectrum of
the Laplacian consisting of such joint eigenfunctions. Furthermore, if
the spectrum is simple, eigenfunctions of the Laplacian are automatically
eigenfunctions of all Hecke operators.

It is no exaggeration to say that these joint eigenfunctions of the
Laplacian and all Hecke operators, usually referred to as Hecke-Maass
cusp forms are of great interest in number theory. Much more is known
about them than about arbitrary eigenfunctions of the Laplacian on a
hyperbolic surface, and so considerable interest has been centered on
the following special case of Rudnick and Sarnak's quantum unique
ergodicity conjecture:

\begin{Question}[Arithmetic Quantum Unique Ergodicity]
Let $M$ be $\Gamma \backslash \H$ for $\Gamma$ a congruence that this
as above, and let $\phi_i$ a sequence of linearly independent joint
eigenfunctions of the Laplacian $\Delta$ and all Hecke operators on
$M$, normalized to have $L^2$-norm one. What are the possible
microlocal lifts of weak$^*$ limits of the probability measures
$\tilde \mu_i$ defined as above by $\tilde \mu_i ( A ) = \int_A
\absolute { \phi_i ( x ) } d \haaron M ( x )$ (we shall refer to these
lifts as arithmetic quantum limits)?
\end{Question}

We review what was known previously regarding this question. T. Watson
\cite{Watson-PhD} has shown that assuming GRH, the only possible
arithmetic quantum limit is the uniform measure on $\Gamma \backslash
\PGL ( 2 , \R )$. In fact, GRH gives an optimal rate of convergence.

Unconditionally, Sarnak has been able to show that a certain
subsequence of the Hecke Maass forms, the CM-forms, satisfies this
conjecture \cite{Sarnak-CM}.

Also unconditionally,
Rudnick and Sarnak \cite{Rudnick-Sarnak} showed that arithmetic quantum
limits cannot be supported on finitely many closed geodesics; similar
ideas were used by the author and Wolpert \cite{Lindenstrauss-imrn,
Wolpert-modulus}.
A small variation of \cite{Lindenstrauss-imrn} gives that any
arithmetic quantum limit is Hecke recurrent for every prime $p$.
Refining these methods, jointly with J. Bourgain, we have been able to prove the following:

\begin{Theorem}[Bourgain and L. \cite{Bourgain-Lindenstrauss}]
\label{theorem with Bourgain}
Every ergodic component of an arithmetic quantum limit has entropy $>
\frac 29$ (in this normalization, the entropy of the uniform Lebesgue
measure is 2).
\end{Theorem}

Using
Theorem~\ref{theorem with Bourgain} together with the Hecke recurrence  we can now deduce from Theorem~\ref{theorem for que} the following theorem regarding arithmetic quantum unique ergodicity

\begin{Theorem}[Arithmetic Quantum Unique Ergodicity]
Let $M$ be $\Gamma \backslash \H$ for $\Gamma$ a congruence lattice as
above. Then any arithmetic quantum limit for $M$ is $c \haaron M$ for
some $c \in [ 0 , 1 ]$; in the compact case, $c = 1$.
\end{Theorem}

This (almost) proves the arithmetic case of the quantum unique
ergodicity conjecture of Rudnick and Sarnak (the only missing piece is
showing $c = 1$ also in the non compact case).

Using the same general strategy, L. Silberman and A. Venkatesh have been able to prove arithmetic quantum unique ergodicity for other $\Gamma \backslash G / K$, specifically for locally symmetric spaces arising from division algebras of prime degree. While the strategy remains the same, several new ideas are needed for this extension, in particular a new micro-local lift for higher rank groups \cite{Silberman-Venkatesh-I}.
\nocite{Ledrappier-Lindenstrauss}

\bibliographystyle{halpha}

\bibliography{newbib}

\end{document}